\documentclass[12pt]{article}
\usepackage{epsfig,psfrag,amsmath,amssymb,latexsym}
\usepackage{amscd }
\usepackage{amsfonts}
\usepackage{graphicx}
\newtheorem{theorem}{Theorem}[section]

\newtheorem{claim}{Claim}[section]

\newtheorem{lemma}{Lemma}[section]

\newenvironment{proof}[1][Proof]{\textbf{#1.} }{\ \rule{0.5em}{0.5em}}

\numberwithin{equation}{section}
\newcommand{\eop}{\mbox{ \vrule height7pt width7pt depth0pt}}
\begin{document}
%------------------
\title{The rate of the convergence of the mean score in random sequence comparison}
\author{J. Lember\thanks{Supported by the Estonian Science Foundation
Grant nr. 7553 and the German Science Foundation (DFG) through CRC 701 at Bielefeld University},   H.
Matzinger, F. Torres\thanks{Supported by the German Science Foundation (DFG) through CRC 701 at Bielefeld University}}

%----------
\maketitle
%-----------
%-------------
\vspace{-0.5cm}
\hspace{0.6cm}\vbox{{\footnotesize
\noindent J\"uri Lember, Tartu University, Institute of Mathematical Statistics\\
Liivi 2-513 50409, Tartu, Estonia. {\it E-mail address:} jyril@ut.ee
\\
Heinrich Matzinger, Georgia Tech, School of Mathematics\\
Atlanta, Georgia 30332-0160, U.S.A. {\it E-mail address:} matzing@math.gatech.edu
\\
Felipe Torres, University of Bielefeld, Faculty of Mathematics \\
Postfach 100131 - 33501 Bielefeld, Germany. {\it E-mail address:} ftorres@math.uni-bielefeld.de
}}
%----------
\vskip 1\baselineskip \noindent \abstract{\noindent We consider a
general class of superadditive scores measuring the similarity of
two independent sequences of $n$ i.i.d. letters from a finite
alphabet. Our object  of interest is the mean score by letter $l_n$.
By subadditivity $l_n$ is nondecreasing and converges to a
limit $l$. We give a simple method of bounding the difference
$l-l_n$ and obtaining the rate of convergence. Our result
generalizes the previous result of Alexander \cite{Alexander1},
where only the special case of the longest common subsequence was
considered.}
%-----------------------------
%--------------
\paragraph{AMS.} 60K35, 41A25, 60C05
%---------------------------
%---------------------------
\section{Introduction}
\begin{center}
\begin{figure}[!h]
\hspace{0.8cm}
\includegraphics[width=15cm]{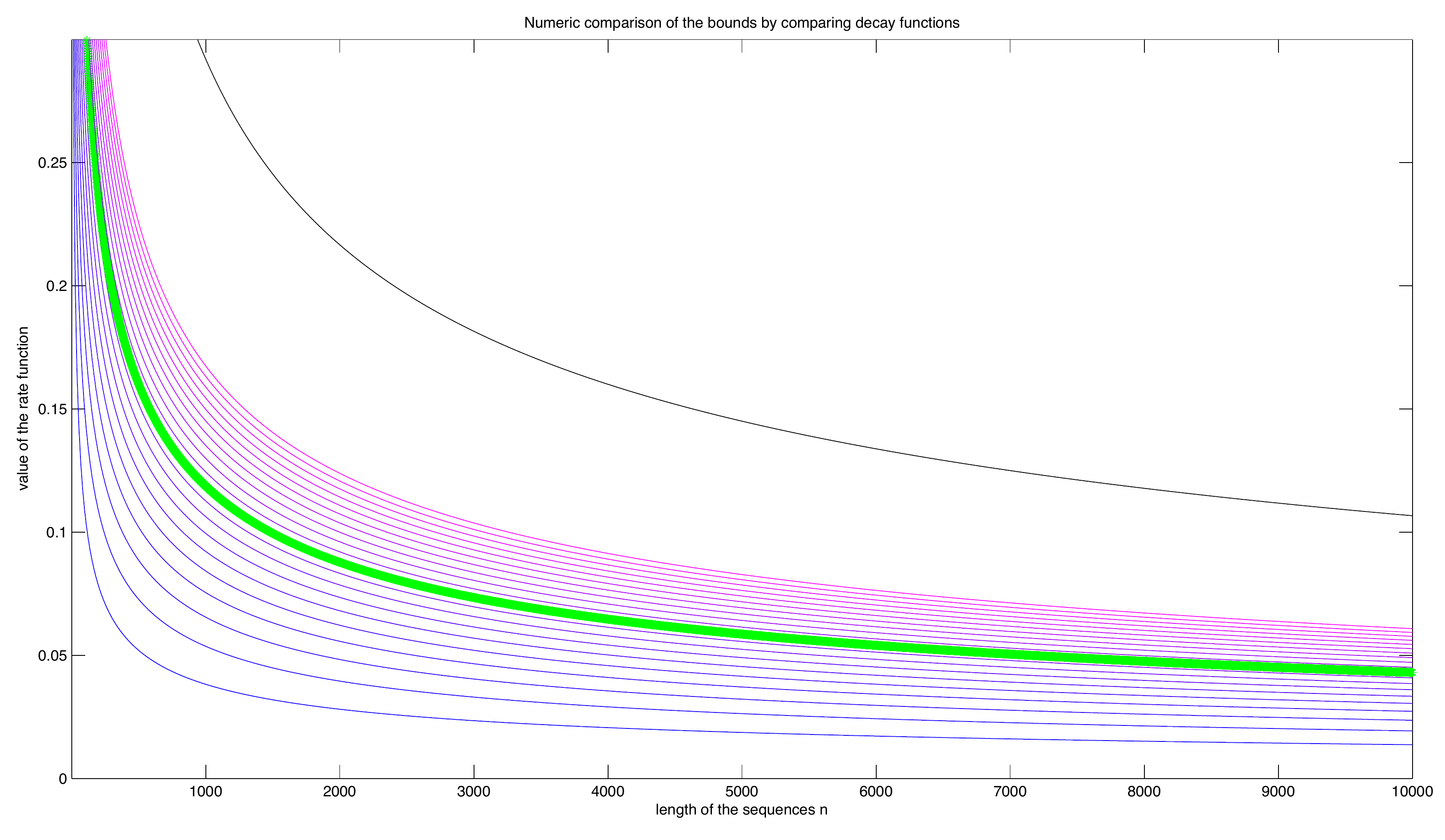}
\caption{\label{fig1} comparison of the bounds (\ref{bestrate}) and (\ref{alexander}) through the functions (\ref{bLMT}) and (\ref{bA}), respectively.}
\end{figure}
\end{center}

Throughout this paper $X_1,X_2,\ldots$ and $Y_1,Y_2,\ldots$ are two
independent sequences of i.i.d. random variables drawn from a finite
alphabet $\mathbb{A}$ and having the same distribution. Since we
mostly study the finite strings of length $n$, let
$X=(X_1,X_2,\ldots X_n)$ and let $Y=(Y_1,Y_2,\ldots Y_n)$ be the
co\-rres\-pon\-ding $n$-dimensional random vectors. We shall usually
refer to $X$ and
$Y$ as random sequences.\\
%-------
The problem of measuring the similarity of $X$ and $Y$ is central in
many areas of applications including computational molecular biology
\cite{christianini, Durbin, Pevzner, SmithWaterman,
watermanintrocompbio} and computational linguistics
\cite{YangLi,LinOch,Melamed1,Melamed2}. In this paper, we consider a
general scoring scheme, where $S:\mathbb{A}\times
\mathbb{A}\rightarrow \mathbb{R}^+$ is a {\it pairwise scoring
function} that  assigns a score to each couple of letters from
$\mathbb{A}$. We assume $S$ to be symmetric and we denote by $F$ and
$A$ the largest  possible score and the largest possible change of
score by one variable, respectively. Formally (recall that $S$ is
symmetric)
$$
F:=\max_{a,b \in \mathbb{A} } S(a,b),\quad A:=\max_{a,b,c \in
\mathbb{A} }|S(a,b)-S(a,c)|.$$ An {\it alignment} is a pair
$(\pi,\mu)$ where $\pi=(\pi_1,\pi_2,\ldots,\pi_k)$ and
$\mu=(\mu_1,\mu_2,\ldots,\mu_k)$ are two increasing sequences of
natural numbers, i.e. $1\leq \pi_1<\pi_2<...<\pi_k\leq n$ and $1\leq
\mu_1<\mu_2<\ldots<\mu_k\leq n.$ The integer $k$ is the number of
aligned letters and $n-k$ is the number of {\it gaps} in the
alignment. Note that our definition of gap slightly differs from the
one that is commonly used in the sequence alignment literature,
where a gap consists of maximal number of consecutive {\it indels}
(insertion and deletion) in one side. Our gap actually corresponds
to a pair of indels, one in $X$-side and another in $Y$ side. Since
we consider the sequences of equal length, to every indel in
$X$-side corresponds an indel in $Y$-side, so considering them
pairwise is justified. In other words, the number of gaps in our
sense is the number of indels in one sequence. We also consider a
{\it gap price} $\delta$. Given the pairwise scoring function $S$ and the
gap price $\delta$, the score of the alignment $(\pi,\mu)$ when
aligning $X$ and $Y$ is defined by
$$U_{(\pi,\mu)}(X,Y):=\sum_{i=1}^kS(X_{\pi_i},Y_{\mu_i})+\delta (n-k).$$
In our general scoring scheme $\delta$ can also be positive,
although usually $\delta\leq 0$ penalizing the mismatch (in this
case $-\delta$ is usually called the gap penalty). We
naturally assume $\delta\leq F$.\\
%-------------
The (optimal)  score of $X$ and $Y$ is defined to be best score over
all possible alignments, i.e.
%$$L_n:=\max_{(\pi,\mu)}U_{(\pi,\mu)}(X,Y).$$
$$L_n:= L(X;Y):=\max_{(\pi,\mu)}U_{(\pi,\mu)}(X,Y).$$
The alignments achieving the maximum are called {\it optimal}. Such
a similarity criterion is most commonly used in sequence comparison
\cite{Watermanphase,Durbin,SmithWaterman,watermanintrocompbio,Vingron}.
When $S(a,b)=1$ for $a=b$ and $S(a,b)=0$ for $a\ne b$, then for
$\delta=0$ the optimal score is equal to the length of the {\it
longest
common subsequence} (LCS) of $X$ and $Y$.\\
It is well-known that  the  sequence $EL_n$, $n=1,2\ldots$ is
superadditive, i.e. $EL_{n+m}\geq EL_n+EL_m$ for all $n,m\geq 1$.
Hence, by Fekete's lemma the ratios $l_n:={EL_n\over n}$ are
nondecreasing and converge to the limit
$$l:=\lim_n l_n=\sup_n l_n.$$
In fact, from Kingman's subadditivity ergodic theorem, it follows
that $l$ is also the a.s. limit of ${L_n\over n}.$ The limit $l$ (which for the LCS-case is called {\it Chvatal-Sankoff constant}) is not
known exactly even for the simplest scoring scheme and the simplest
model for $X$ and $Y$, so it is usually estimated by simulations. Using
McDiarmid's inequality (see (\ref{McD2})) one can estimate $l_n$ with prescribed
accuracy; to obtain confidence
intervals for $l$, the difference $l-l_n$ should be
estimated. This is the aim of the present paper. \\
To our best knowledge, the difference $l-l_n$ has been theoretically studied only by Alexander in \cite{Alexander1}, though there exist many numeric results on the value of $l_n$ or its distribution in various contexts \cite{Baeza1999,BMNWilson,Sankoff1,Paterson1,Deken,FuLuo,meancurve,kiwi,Paterson2}. Alexander proved that in the case of the LCS, for any $C>(2+\sqrt{2})$ there exists an integer $n_o(C)$ such that
\begin{equation}\label{alexander}
l-l_n\leq C\sqrt{{\log n \over n}},\quad \text{provided   }
n>n_o(C).
\end{equation}
The bound (\ref{alexander}) is independent of the common law of $X$
and $Y$, and the integer $n_o(C)$ can be exactly determined. Hence
the bound (\ref{alexander}) can be used for the calculation of
explicit confidence intervals.\\\\
%-----------------
Our main result is the following:
%------------
\begin{theorem}\label{thm} Let $n\in \mathbb{N}$ be even. Then, with any $c>\sqrt{A}$,
\begin{equation}\label{bestrate}
l-l_{n}\leq c\sqrt{{2\over n-1}\Big( {n+1\over n-1}+\ln (n-1)
\Big)}+{F\over n-1}.\end{equation}
\end{theorem}
Note that by the monotonicity of  $l_n$, the  assumption on $n$ even 
actually is not restrictive. In fact, Alexander's main result
(Prop. 2.4 in \cite{Alexander1}) is also proven for $n$ even.
Theorem \ref{thm} and its proof generalize Alexander's result in
many ways:
\begin{enumerate}
\item Theorem \ref{thm} applies for a general scoring scheme, not
just for the LCS. This is due to the fact that our proof is based
solely on McDiarmid's large deviation equality, whilst Alexander's
proof, although using also  McDiarmid's inequality, is mainly based on first passage percolation techniques.
Despite the fact that the percolation approach  applies in many other
situations rather than sequence comparison (see \cite{Alexander2}), it is not clear
whether it can be efficiently applied to our general scoring
scheme. For McDiarmid's inequality, however, it makes no
difference what kind of scoring is used. This gives us
reasons to believe that our proof is somehow "easier" than the one in
\cite{Alexander1}.
 %------------------
\item The proof of Theorem \ref{thm} relates the rate of the
convergence of $l_n$ to  the cardinality of the set of partitions ${\cal
B}_{k,n}$ (see Lemma \ref{lemma}) so that finding the good rate boils down to the good estimation of
$|{\cal B}_{k,n}|$. The bound (\ref{bestrate}) corresponds to a
particular estimate of $|{\cal B}_{k,n}|$, any better estimate would give a sharper bound
and, probably, also a faster rate. In a sense, the cardinality
$|{\cal B}_{k,n}|$ could be interpreted as the complexity of the
model and  the  relation between the rate of convergence and the
complexity of the model is a well-known fact in statistics (see e.g.
\cite{BarronBirgeMassart}).
\item \label{LCSbetter} When applied to the LCS, our bound (\ref{bestrate}) is
sharper than (\ref{alexander}). Indeed, for the case of LCS the
constants $A$ and $F$ in (\ref{bestrate}) can be taken equal to
one and the smaller  constants make the difference.  In other words, for the case of LCS  both results yield the
rate $C\sqrt{{\ln n\over n}}$, but the constant $C$ is different
($C>3.42$ in Alexander's result and $c>\sqrt{2}$ in ours).
\end{enumerate}
%-----------------
We can easily compare (\ref{bestrate}) and (\ref{alexander}) by comparing the decay of the two following functions:
\begin{eqnarray}
R : \{1,\dots,10000\} &\to& \mathbb{R}^+ \nonumber \\
R(n) &=&  (3.42+0.1)\sqrt{\frac{\ln n}{n}}\label{bA} \\
Q_F : \{ 1,\dots,10000\} \times \{ 0.1,\dots,2\} &\to& \mathbb{R}^+ \nonumber \\
Q_F(n,A) &=&  \sqrt{A}\,\, \sqrt{{2\over n-1}\Big( {n+1\over n-1}+\ln (n-1)\Big)}+{F\over n-1} \label{bLMT} \nonumber \\
\end{eqnarray}
In figure \ref{fig1}, we can see the improved bound (\ref{bestrate}) given by function (\ref{bLMT}) (changes of $A$ are represented in colours, $F=1$) over the bound by Alexander (\ref{alexander}) given by function (\ref{bA}) (in black). Note that the dark blue curve corresponds to $A=0.1$ whilst the light violet curve to $A=2$ (namely, the colour gets lighter as $A$ increases). The curve in green corresponds to the case $A=1$ (ie, our bound for the LCS case).
%-----------------
\section{Confidence bounds for $l$}
Suppose that $k$ samples of $X^i=X_1^i,\dots,X_n^i$ and
$Y^i=Y_1^i,\dots,Y_n^i$, $i=1,\dots,N$ are generated. Let $L^i_n$ be
the score of the $i$-th sample. Thus $EL_n^i=n \,l_n$. By McDiarmid's
inequality (see (\ref{McD}) below), for every $\rho >0$
\begin{equation}
P\left(\frac{1}{kn} \sum_{i=1}^k L^i_n - l_n < -\rho\right) =
P\left(\sum_{i=1}^k L^i_n - kn\,l_n < -kn\rho\right) \le \exp\left[
-\frac{\rho^2kn}{A^2}\right].
\end{equation}
Let
\[ \bar{L}_n:=\frac{1}{kn} \sum_{i=1}^k L^i_n.\]
If $n$ is even, by (\ref{bestrate}) and (\ref{bLMT}) we have that $l \le l_n + Q_F(n,A)$ and then
\begin{equation}
P(\bar{L}_n + \rho + Q_F(n,A) \ge l) \ge P(\bar{L}_n + \rho \ge l_n)
= P(\bar{L}_n - l_n \ge -\rho) \ge 1-\exp\left[ -\frac{\rho^2
kn}{A^2}\right].
\end{equation}
Now, given $\varepsilon >0$, choose $\rho=\rho(n,\varepsilon)$ so that the right hand side in the last inequality is equal to $1-\varepsilon$:
\[ \rho(n,\varepsilon) = A \sqrt{\frac{\ln(1/\varepsilon)}{kn}}.\]
So, with probability $1-\varepsilon$, we obtain one side confidence interval as follows:
\begin{equation}\label{PAC1}
l \le \bar{L}_n + Q_F(n,A) + A\sqrt{\frac{\ln(1/\varepsilon)}{kn}}.
\end{equation}
In statistical learning, the inequalities of type (\ref{PAC1}) are known as PAC inequality (probably almost correct inequalities).
The two-sided confidence bounds are, with probability $1-\varepsilon$, as follows:
\begin{equation}\label{PAC2}
\bar{L}_n -A\sqrt{\frac{\ln(2/\varepsilon)}{kn}} \le l \le \bar{L}_n
+ Q_F(n,A) + A\sqrt{\frac{\ln(2/\varepsilon)}{kn}}.
\end{equation}
The bounds in (\ref{PAC2}) suggest to use the estimate
\[ \hat{l}_n := \bar{L}_n+\frac{Q_F(n,A)}{2}\]
so that the confidence bounds for this estimate are
\begin{equation}\label{confidencebounds}
P\left( |\hat{l}_n - l| \le A\sqrt{\frac{\ln(2/\varepsilon)}{kn}} +
\frac{Q_F(n,A)}{2}\right) \ge 1-\varepsilon.
\end{equation}
Alexander \cite{Alexander1} obtained, for $n=100000$, $k=2$ and $A=F=1$ (for the LCS case), the following bounds:
\begin{equation}\label{alexandercbounds}
P(| \hat{l}_n - l| \le 0.0264) \ge 0.95.
\end{equation}
By using (\ref{confidencebounds}) and (\ref{bLMT}) we obtain, for $n=100000$, $k=2$ and $A=F=1$ (for the LCS case), the following bounds:
\begin{equation}\label{LMTcbounds}
P(| \hat{l}_n - l| \le 0.0122) \ge 0.95.
\end{equation}
%{\bf Felipe: you had forgotten 2 in two-sided confidence bounds.
%That is why rthe original number 0.0119 is not correct. I roughly
%corrected it, but check the calculations!}\\\\
%-----------------------
It is clear that (\ref{LMTcbounds}) is sharper than (\ref{alexandercbounds}). To our best knowledge, the best previous lower and upper bounds for $l$, in the LCS context for $\mathbb{A}=\{0,1\}$, were due to Dancik \cite{Dancik}, Dancik and Paterson \cite{Paterson1,Paterson2} (0.773911 and 0.837623, respectivley) and Lueker \cite{Lueker} (0.788071 and 0.826280, respectively).
\\
\\
{\bf Remark: }The inequality (\ref{PAC1}) confirms the well-known fact that it is better to generate one sample of length $kn$ rather than $k$ samples of length $n$. Indeed, with one sample of length $kn$, the inequality (\ref{PAC1}) becames
\begin{equation}
l \le \bar{L}_n + Q_F(kn,A)+ A\sqrt{\frac{\ln(1/\varepsilon)}{kn}}
\end{equation}
and since $Q_F(kn,A) < Q_F(n,A)$, the bounds get narrower.
%-----------------
\section{Proof of the main result}
\subsection{The set of partitions ${\cal B}_{k,n}$}
In this section, we shall consider the sequences $X$ and $Y$ with
length $kn$ where $k,n$ are nonnegative integers. Let $(\pi,\mu)$ be
an arbitrary alignment of $X$ and $Y$. Let
$\nu=(\nu_1,\ldots,\nu_{r+1})$ and $\tau=(\tau_1,\ldots,\tau_{r+1})$
be  vectors satisfying
\begin{equation}\label{part1}
1=\nu_1\leq \nu_2 \leq \ldots \leq \nu_r\leq \nu_{r+1}=kn+1,\quad 1=
\tau_1\leq \tau_2\leq \ldots \leq \tau_r\leq \tau_{r+1}=kn+1.
\end{equation}
We say that
the pair $(\nu,\tau)$ forms a {\bf $r$-partition of the alignment
$(\pi,\mu)$}
if for any $j=1,\ldots,r$, the following conditions are simultaneously satisfied:\\
1)  if, for some $i=1,\ldots k$, it holds that $\nu_j \leq \pi_i < \nu_{j+1}$, then $\tau_j \leq \mu_i < \tau_{j+1}$;\\
2)  if, for some $i=1,\ldots k$, it holds that $\tau_j \leq \mu_i < \tau_{j+1}$, then $\nu_j \leq \pi_i < \nu_{j+1}$.\\\\
Thus
$(\nu,\tau)$ is a $r$-partition, if the sequences $X$ and $Y$ can be partitioned into $r$ pieces
\begin{align*}
&(X_1,\ldots,X_{\nu_2-1}),(X_{\nu_2},\ldots,X_{\nu_3-1}),\ldots,(X_{\nu_r},\ldots,X_{kn})\\
&(Y_1,\ldots,Y_{\tau_2-1}),(Y_{\tau_2},\ldots,Y_{\tau_3-1}),\ldots,(Y_{\tau_r},\ldots,Y_{kn})
\end{align*}
such that the alignment  $(\pi,\mu)$ aligns a piece
$(X_{\nu_j},\ldots,X_{\nu_{j+1}-1})$ with the piece
$(Y_{\tau_j},\ldots,Y_{\tau_{j+1}-1})$, where $j=1,\ldots r$. It is
important to note that the pieces might be empty, i.e. it might be
that  $\nu_j=\nu_{j+1}$ (or $\tau_j=\tau_{j+1}$), meaning that
$({\tau_j},\ldots,{\tau_{j+1}-1})$ cannot contain any elements of
$\mu$, otherwise the requirement 2) would be violated (or
$({\mu_j},\ldots,{\mu_{j+1}-1})$ cannot contain any elements of
$\tau$, otherwise the requirement 1) would be violated). Hence, if
for a partition a piece  of $X$ is empty, then the corresponding
piece of $Y$ cannot have any aligned letter.\\\\
%-----------------
The following observation shows that any alignment of $X$ and $Y$
can be partitioned into $r$ pieces such that $k\leq r\leq
\lceil{2kn\over 2n-1}\rceil$ and such that the sum of the lengths of
aligned pairs in each partition is always at most $2n$.
We believe that the idea of the proof as well as the meaning of the partition becomes transparent by an example.\\\\
%------------------
{\bf Example.} Let $n=3,k=4$. Let $\pi =(1,5,6,9,10,12)$ and $\mu
=(2,3,4,6,9,10)$. The alignment $(\pi,\mu)$ can be represented  as
follows
\begin{center}
% use packages: array
\begin{tabular}{l|l|l|l|l|l|l|l|l|l|l|l|l|l|l|l|l|l|l}
$X$ &- & 1 & 2 & 3 & 4 & 5 & 6 & 7 & 8 & - & 9 & - & - & 10 & 11 & 12 & - & - \\
\hline $Y$ & 1 & 2 & - & - & - & 3 & 4 & - & - & 5 & 6 & 7 & 8 & 9 &
- & 10 & 11 & 12
     \end{tabular}
     \end{center}
The table above indicates that $X_1$ is aligned with $Y_2$, $X_5$ is
aligned with $Y_3$ and so on; the rest of the letters
are unaligned, so we say that they are aligned with gaps.  In the
table, there are two types of columns: the columns with two figures
(aligned pairs) and the columns with one figure (unaligned pairs).
Let $u_i\in \{1,2\}$ be the number of figures in the $i$-th column,
and let $s_j=u_1+\cdots + u_j$ be the corresponding cumulative sum.
To get an $r$-partition proceed as follows: start from the
beginning of the table (most left position) and find $j$ such that $s_j=2n$. Since the
cumulative sum increases by one or two, such a $j$ might not
exist. In this case find $j$ such that $s_j=2n-1$. In the present
example $n=3$, thus we are looking for $j$ such that $s_j=6$. Such a
$j$ is 5. The first five columns thus form the first part of the
partition and there are exactly $2n=6$ elements in the first part (those elements are $X_1,X_2,X_3,X_4,Y_1$ and $Y_2$).
Now disregard  the first five columns from the table and start the
same procedure afresh. Then the second part is obtained and so on.
In the following table the vertical lines indicate the different
parts obtained by the aforementioned procedure: the first two parts have
six elements, the third and fourth  has five elements and the last
part consists of one element:
\begin{center}
% use packages: array
\begin{tabular}{l| l l l l l |l l l l |l l l l |l l l l |l}
$X$ &- & 1 & 2 & 3 & 4  & 5 & 6 & 7 & 8 & - & 9 & - & - & 10 & 11 & 12 & - & - \\
\hline
$Y$ & 1 & 2 & - & - & - & 3 & 4 & - & - & 5 & 6 & 7 & 8 & 9 &
- & 10 & 11 & 12
     \end{tabular}
     \end{center}
From the table, we read the corresponding pieces from  the $X$-side:
$(1,4),(5,8),(9,9),$ $(10,12), \emptyset$ as well as the ones from
the $Y$-side: $(1,2),(3,4),(5,8),(9,11),(12,12)$. The
co\-rres\-pon\-ding vectors $\nu$ and $\tau$  are thus
$\nu=(1,5,9,10,13,13)$, $\tau=(1,3,5,9,12,13)$. The number of parts
in such a partition is clearly at least $k$ (corresponding to the
case that all pairs sum up to $2n$) and at most $\lceil {2kn\over
2n-1}\rceil$ (corresponding to the case that all pairs except the
last one sum up to $2n-1$). In our example is $r=5=\lceil{24\over
5}\rceil$. Now, it is clear that the following claim holds.
%----------------------
\begin{claim}\label{claim5}
Let $X$, $Y$ be sequences of length $kn$ and let $(\pi,\mu)$ be an
arbitrary alignment of $X$ and $Y$. Then there exist an integer $r$ such that
$k\leq r \leq \lceil{2kn\over 2n-1}\rceil$ and an
$r$-partition $(\nu,\tau)$ of $(\pi,\mu)$ such that for every
$j=1,\ldots,r-1$, it holds
\begin{equation}\label{Pcond5}
(\nu_{j+1}-\nu_j)+(\tau_{j+1}-\tau_j)\in \{2n, 2n-1\}\quad
\text{and}\quad (\nu_{r+1}-\nu_r)+(\tau_{r+1}-\tau_r)\leq 2n.
\end{equation}
\end{claim}
%-------------
Let, for every $r$, ${\cal B}^r_{k,n}$ be the set of vectors
$\nu=(\nu_1,\ldots,\nu_{r+1})$ and $\tau=(\tau_1,\ldots,\tau_{r+1})$
satisfying (\ref{part1}) and (\ref{Pcond5}). Let
$${\cal B}_{k,n}=\bigcup_{r=k}^{\lceil{2kn\over 2n-1}\rceil}{\cal B}^r_{k,n}.$$
We shall call the elements of ${\cal B}_{k,n}$ as the partitions.
For every partition
$(\nu,\tau)\in {\cal B}^r_{k,n}$, we define
$$L_{kn}(\nu,\tau):=\sum_{i=1}^r
L(X_{\nu_j},\ldots,X_{\nu_{j+1}-1};Y_{\tau_j},\ldots,Y_{\tau_{j+1}-1}),$$
where
$L(X_{\nu_j},\ldots,X_{\nu_{j+1}-1};Y_{\tau_j},\ldots,X_{\tau_{j+1}-1})$
is the optimal score between $X_{\nu_j},\ldots,$ $X_{\nu_{j+1}-1}$
and $Y_{\tau_j},\ldots,Y_{\tau_{j+1}-1}$. The key observation is the
following: if $(\pi,\mu)$ is optimal for $X,Y$ and $(\nu,\tau)$ is a
$r$-partition of $(\pi,\mu)$, then $L_{kn}=L_{kn}(\nu,\tau).$ By Claim
\ref{claim5}, every alignment, including the optimal one, has at
least one partition from the set ${\cal B}_{k,n}$, hence  it follows
that
\begin{equation}\label{max}
L_{kn}=\max_{(\nu,\tau)\in {\cal B}_{k,n}} L_{kn}(\nu,\tau).
\end{equation}
%------------
\begin{claim}\label{claim6}
For every $r$-partition $(\nu,\tau)\in {\cal B}_{k,n}$,
\begin{equation}\label{E6}
E\big(L_{kn}(\nu,\tau)\big)\leq {r\over 2} EL_{2n}\leq {1\over 2}\Big\lceil{2kn\over 2n - 1}\Big\rceil EL_{2n}.
\end{equation}
\end{claim}
\begin{proof}
Let $(\nu,\tau)\in {\cal B}^r_{k,n}$ with $r \leq \lceil{2nk\over 2n-1}\rceil$. Let $j$ be such that
$(\nu_{j+1}-\nu_j)+(\tau_{j+1}-\tau_j)=2n$. Thus, there exists an
integer $u\in \{-n,\ldots, n\}$ such that $\nu_{j+1}-\nu_j=n-u$ and
$\tau_{j+1}-\tau_j=n+u$.  Since $X_1,X_2,\ldots, Y_1,Y_2,\ldots$ are
i.i.d., we have
\begin{align*}
&E\big(L(X_{\nu_j},\ldots,X_{\nu_{j+1}-1};Y_{\tau_j},\ldots,Y_{\tau_{j+1}-1})\big)=E\big(L(X_1,\ldots,X_{n-u};Y_1,\ldots,Y_{n+u})\big)=\\
&E\big(L(X_{n-u+1},\ldots,X_{2n};Y_{n+u+1},\ldots,Y_{2n})\big)\leq
{1\over
2}E\big(L(X_{1},\ldots,X_{2n};Y_{1},\ldots,Y_{2n})\big)={1\over 2}
EL_{2n}.
\end{align*}
The last inequality follows from the superadditivity:
\begin{eqnarray}
L(X_1,\ldots,X_{n-u};Y_1,\ldots,Y_{n+u})+L(X_{n-u+1},\ldots,X_{2n};Y_{n+u+1},\ldots,Y_{2n})&& \nonumber \\
\leq L(X_{1},\ldots,X_{2n};Y_{1},\ldots,Y_{2n}). && \nonumber
\end{eqnarray}
If $(\nu_{j+1}-\nu_j)+(\tau_{j+1}-\tau_j)<2n$, then by the same
argument
\begin{align*}
&E\big(L(X_{\nu_j},\ldots,X_{\nu_{j+1}-1};Y_{\tau_j},\ldots,Y_{\tau_{j+1}-1})\big)\leq
E\big(L(X_1,\ldots,X_{n-u};Y_1,\ldots,Y_{n+u})\big)\leq {1\over 2}
EL_{2n}.
\end{align*}
Hence the first inequality in (\ref{E6}) follows. The second
inequality follows from the condition $r \leq \lceil{2nk\over 2n-1}\rceil$.
%definition of ${\cal B}_{k,n}$.
\end{proof}
%------------------
\subsection{The size of ${\cal B}_{k,n}$ and the rate of
convergence}
%-----------
In the following we prove the main theoretical result that links the
rate of the convergence to the rate at which the number of elements
in $|{\cal B}_{k,n}|$ grows as $k$ increases. Our proof is entirely
based on McDiarmid's inequality, so let us recall it for the sake of
completeness: Let $Z_1,\ldots,Z_{2m}$ be independent random
variables and $f(Z_1,\ldots,Z_{2m})$ be a function so that changing
one variable changes the value at most $A$. Then for any $\Delta>0$,
\begin{equation}\label{McD}
P\Big(f(Z_1,\ldots,Z_{2m})-Ef(Z_1,\ldots,Z_{2m})>\Delta\Big) \leq
\exp\left[-{{\Delta}^2\over m A^2}\right].
\end{equation}
For the proof, we refer \cite{lugosi}. We apply (\ref{McD}) with $L$
in the role of $f$ to the independent (but not necessarily
identically distributed) random variables
$X_1,\ldots,X_{m},Y_1,\ldots,Y_{m}$. It is easy but important to see
that independently of the value of $\delta$, changing one random
variable changes the score at most by $A$ so that in our case
(\ref{McD}) is
\begin{equation}\label{McD2}
P\Big(L_m-EL_m>\Delta\Big) \leq \exp\left[-{{\Delta}^2\over m
A^2}\right].
\end{equation}
%------------------------
%-----------------
\begin{lemma}\label{lemma}
Suppose that for any  $n$ and for $k$ big enough
\begin{equation}\label{assumption2}
|{\cal B}_{k,n}|\leq \exp[\big(\psi(n)+o(k)\big)kn],
\end{equation}
where $\psi(n)$ does not depend on $k$. Let $u(n)>A\sqrt{\psi(n)}$.
Then
\begin{equation}\label{statement}
l-l_{2n}\leq u(n)+{l_{2n}\over 2n-1}\leq u(n)+{l\over 2n-1}\leq
u(n)+{F\over 2n-1}.\end{equation}
\end{lemma}
%-----------------------
\begin{proof}
%----------------------------------
Let $(\nu,\tau)\in {\cal B}_{k,n}$. Recall (\ref{E6}). Thus, from
(\ref{McD2}), we get that for any $\rho>0$,
\begin{equation}\label{McD4}
P\Big(L_{kn}(\nu,\tau)-{1\over 2}\left\lceil {2kn\over 2n-1}\right\rceil EL_{2n} > \rho kn \Big)\leq P\Big(L_{kn}(\nu,\tau)- E\big(L_{kn}(\nu,\tau)\big)
\rho kn \Big)\leq \exp\left[-{\rho^2kn\over A^2}\right].
\end{equation}
From (\ref{max}) and  (\ref{assumption2})  it now follows that, for
big $k$
\begin{align*}
P\left({L_{kn}\over kn}-{1\over k}\left\lceil {2kn\over 2n-1}\right\rceil  l_{2n} > \rho \right)&\leq \sum_{(\nu,\tau)\in {\cal
B}_{k,n}}P\Big(L_{kn}(\nu,\tau)-{1\over 2} \left\lceil {2kn\over 2n-1}\right\rceil  EL_{2n} >\rho kn \Big)\\
&\leq |{\cal B}_{k,n}|\exp\left[-{\rho^2kn\over A^2}\right]\leq
\exp\Big[\Big(\psi(n)+o(k)-\left({\rho\over A}\right)^2\Big)kn\Big].
\end{align*}
We consider $n$ fixed and let $k$ go to infinity. If
$u(n)>A\sqrt{\psi(n)}$, then there exists $K(n)<\infty$ so that for
every $k>K(n)$,
\[\psi(n)+o(k)-\left({u(n)\over A}\right)^2<\,\,{1\over 2}\Big(\psi(n)-\left({u(n)\over A}\right)^2\Big).\]
Hence, replacing in the inequalities above $\rho$ with
$u(n)$, we obtain for every  $k>K(n)$,
\begin{equation}\label{main2}
P\left({L_{kn}\over kn}-{1\over k}\left\lceil {2kn\over 2n-1}\right\rceil l_{2n} > u(n)\right)\leq \exp\left[{1\over 2}\left(\psi(n)-\left({u(n)\over
A}\right)^2\right)nk\right]=\exp[-d_nk],
\end{equation}
where
$$d_n:=\left( \left({u(n)\over A}\right)^2-\psi(n)\right)n>0.$$
Now recall the assumption that $\delta\leq F$. Hence
for any $n$ and $k$, the random variable ${L_{kn}\over kn}$ is bounded by $F$. From (\ref{main2}), it thus follows that for any
$k$
$$E\Big({L_{kn}\over kn}\Big)=l_{kn}\leq {1\over  k} \left\lceil {2kn\over 2n-1} \right\rceil l_{2n} +u(n)+F\exp[-d_nk].$$
Since $l_{kn}\to l$ as  $k\to \infty$ and
$${1\over  k} \left\lceil {2kn\over 2n-1} \right\rceil\leq {2n\over
2n-1}+{1\over k},$$ we obtain that for any $n$,
$$l\leq \left({2n\over 2n-1}\right)l_{2n}+u(n)=l_{2n}\left(1+{1\over 2n-1}\right)+u(n).$$
\end{proof}
%------------
%-----------------------------
\paragraph{The proof of Theorem \ref{thm}.}
From Lemma \ref{lemma}, it follows that to obtain a bound to
$l-l_n$, a suitable estimator of $|{\cal B}_{k,n}|$ satisfying
(\ref{assumption2}) should
be found. \\
%-----------
Let us estimate $|{\cal B}^r_{k,n}|$. The number of parts in the $X$
side is bounded above by the number of combination with repetition
from $nk+1$ by $r-1$. The repetitions allow empty parts. When the
size of a part in $X$-side is  $m$, then, except from the last part,
the size of the co\-rres\-pon\-ding part on $Y$ side has two
possibilities: $2n-1-m$ or $2n-m$. Hence to any $r$-partition of
$X$-size corresponds at most $2^{r-1}2n$ options in $Y$ side. In the
following we use the fact that the number of combination with
repetition from $nk+1$ by $r-1$ is ${nk+r-1 \choose r-1}$ and for
any non-negative integers $a>b$ it holds $${a\choose b}\leq
\exp\left[h_e\left({b\over a}\right)a\right],$$ where $h_e(q):=-q\ln q-(1-q)\ln(1-q)$ is
the binary entropy function. Since $r\leq \lceil{2nk\over
2n-1}\rceil$ implies that $r-1\leq {2nk\over 2n-1}$, we thus have
for $n\geq 2$
\begin{align*}
|{\cal B}^r_{k,n}|&\leq (2^{r-1}2n){nk+r-1 \choose r-1} \\
&\leq \exp\Big[(r-1)(\ln 2) +\ln (2n)+h_e\left({r-1\over nk+r-1}\right)(nk+r-1)\Big]\\
&\leq \exp\left[\left({\ln 4\over 2n-1}+{\ln (2n)\over nk}+ h_e\left({r-1\over nk+r-1}\right)\left(1+{2\over 2n-1}\right)\right)nk\right]\\
&\leq \exp\left[\left({\ln 4\over 2n-1}+{\ln (2n)\over nk}+h_e\left({2\over 2n+1}\right)\left({2n+1\over 2n-1}\right)\right)nk\right].
\end{align*}
The last inequality follows from the inequalities
$${r-1\over nk+r-1}\leq {{2nk\over 2n-1}\over nk+{2nk\over
2n-1}}={2\over 2n+1}$$ so that if $n\geq 2$, then ${2\over 2n+1}\leq
0.5$ and
$$h_e\left({r-1\over nk+r-1}\right)\leq h_e\left({2\over 2n+1}\right).$$
Hence
\begin{align*}
|{\cal B}_{k,n}|&\leq \left({2nk\over 2n-1}-k+2\right)\exp\left[\left({\ln 4\over
2n-1}+{\ln (2n)\over nk}+h_e\left({2\over 2n+1}\right)\left({2n+1\over 2n-1}\right)\right)nk\right]\\
&=\left({k\over 2n-1}+2\right)\exp\left[\left({\ln 4\over 2n-1}+{\ln (2n)\over nk}+h_e\left({2\over 2n+1}\right)\left({2n+1\over 2n-1}\right)\right)nk\right]\\
&=\exp\left[\ln\left({k\over 2n-1}+2\right)+\left({\ln 4\over 2n-1}+{\ln (2n)\over nk}+h_e\left({2\over 2n+1}\right)\left({2n+1\over 2n-1}\right)\right)nk\right]\\
&=\exp\left[\left({\ln\left({k\over 2n-1}+2\right)+\ln(2n)\over nk}+{\ln 4\over 2n-1}+h_e\left({2\over 2n+1}\right)\left({2n+1\over 2n-1}\right)\right)nk\right]
\end{align*}
\begin{align*}
&=\exp\left[\left(o(k)+{\ln 4\over 2n-1}+h_e\left({2\over 2n+1}\right)\left({2n+1\over 2n-1}\right)\right)nk\right] \\
& \leq \exp\left[\left(o(k)+{2\over 2n-1}\left({2n+1\over 2n-1}+\ln(2n-1)\right)\right)nk\right],
\end{align*}
where the last inequality follows from the inequality
\begin{equation}
h_e\left({2\over 2n+1}\right)\leq {2\over 2n+1}\left({2n+1\over 2n-1}+\ln \left({2n-1\over 2}\right)\right).
\end{equation}
Hence (\ref{assumption2}) holds with
$$\psi(n)={2\over 2n-1}\left({2n+1\over 2n-1}+\ln(2n-1)\right).$$
%$$\psi(n)={1\over 2n-1}\left(h_e\left({2\over 2n+1}\right)\left(2n+1\right)\right)+{\ln 4\over 2n-1}.$$
%----------------------
%With
% we get
%---------------------------------
The inequality (\ref{bestrate}) now follows from Lemma \ref{lemma}$.\quad\eop$

\vspace{36pt}
\noindent {\bf \large Acknowledgments}
\vspace{12pt}

\noindent The authors would like to thank the support of the German Science Foundation (DFG) through the Collavorative Research Center 701 "Spectral Structures and Topological Methods in Mathematics" (CRC 701) at Bielefeld University and the support of Barbara Gentz with the research stay of J. Lember at the CRC 701. Additionally, F. Torres would like to thank the partial support of the International Graduate College "Stochastics and Real World Models" (IRTG 1132) at Bielefeld University.
\newpage
\bibliographystyle{plain}
\bibliography{bio}
\end{document}